\documentclass[12pt]{article}
\title{Equimultiple locus of embedded algebroid surfaces and blowing--up in
arbitrary characteristic}
\author{
Piedra, R.\footnote{Supported by FQM 304 and BFM 2000--1523.}
\\ Depto. de \'Algebra  \\ Universidad de Sevilla \and
Tornero, J.M.\footnote{Supported by FQM 218 and MTM2004--07203--C02.} 
\\ Depto. de \'Algebra  \\ Universidad de Sevilla}
\date{November, 2005}

\newcommand{\NN}{{\bf N}}

\newcommand{\PP}{{\bf P}}

\newcommand{\vs}{\vspace{.15cm}}

\newcommand{\df}{\noindent {\bf Definition.-- }}

\newcommand{\thr}{\noindent {\bf Theorem.-- }}
\newcommand{\lema}{\noindent {\bf Lemma.-- }}
\newcommand{\dem}{{\bf Proof.-- }}
\newcommand{\obs}{\noindent {\bf Remark.-- }}

\newcommand{\ntt}{\noindent {\bf Notation.-- }}

\newcommand{\ord}{\mbox{ord}}
\newcommand{\Sp}{\mbox{Spec}}

\newcommand{\ol}{\overline}

\newcommand{\cE}{{\cal E}}
\newcommand{\cS}{{\cal S}}

\begin{document}

\maketitle

\abstract{This paper extends previous results of the authors,
concerning the behaviour of the equimultiple locus of algebroid
surfaces under blowing--up, to arbitrary characteristic.}

\vs


\section{Introduction}

Let $K$ be an algebraically closed field of characteristic $p$ and
$\cS = \Sp (K[[X,Y,Z]]/(F))$ an embedded algebroid surface. With no loss 
of generality $\cS$ can be considered to be defined
by a Weierstrass equation
$$
F(Z) = Z^n + \sum_{k=0}^{n-1} a_k (X,Y) Z^k,
$$
where $n$ is the multiplicity of $\cS$, that is, $\ord \left(a_k
\right) \geq n-k$ for all $k=0,...,n-1$. After a change of
variables we can assume that either $\ol{F}$, the initial form of
$F$, is not a power of a linear form, or $\ol{F} = Z^n$. From now
on, by a Weierstrass equation we will mean an equation like this.

In this situation the equimultiple locus of $\cS$ is
$$
\cE (\cS)
= \left\{ P \in \Sp(K[[X,Y,Z]]) \; | \; F \in P^{(n)} \right\},
$$
which is never empty, as $M=(X,Y,Z)$ always lies in $\cE (\cS)$.
Note that all elements of $\cE (\cS)$, different from $M$, can 
be assumed to have the form $P=(Z+G(X,Y), H(X,Y))$.

Geometrically speaking, the equimultiple locus represents points
at which the multiplicity is the same than in the origin; hence
they are the ``closest'' points to the origin, in (coarse) terms
of singularity complexity. We will note by $\cE_0 (\cS)$ the
subset of smooth elements of $\cE (\cS)$.

In our previous paper \cite{PT} we prove a theorem relating $\cE_0
\left( \cS^{(1)} \right)$ to $\cE (\cS)$, where $\cS^{(1)}$ is
the blowing--up of $\cS$ centered in an element of $\cE_0 (\cS)$.
Our aim is extending this result to the arbitrary characteristic 
case.

The main difficulty in the positive characteristic case comes from
the fact that, if $p$ divides $n$, we cannot apply the
Tchirnhausen transformation. This deceivingly naive procedure
assures us, when it can be carried out, that:

\begin{enumerate}
\item[(a)] If $\ol{F}$ is the power of a linear form (that is,
if the tangent cone is a plane), then it must be $\ol{F} = Z^n$.
\item[(b)] All the elements of $\cE (\cS)$ contain $Z$.
\item[(c)] If, after blowing--up, the multiplicity remains the
same, the above conditions still hold.
\end{enumerate}

In our current situation, as we pointed out, one finds a regular
parameter verifying (a) with a simple change of coordinates.
Moreover, after Mulay's work (\cite{Mulay}) one can also find a
regular parameter satisfying (b). We should note here that the
outstanding result of Mulay is not constructive (much in the
spirit of Abhyankar's beautiful theory on good points \cite{GP}).
Even though one can manage to find a parameter verifying (a) and
(b), the preservation of (a) and (b) under certain blowing--ups is
false in positive characterisitic, as is well--known (see \cite{RP}
for how this might be overcome in a Levi--Zariski resolution).

The existence of such an amenable parameter would have made
the null characteristic proof valid for $p>0$. The peculiarities of 
the case $p|n$ have forced us to use completely different 
strategies for most parts of the result, although our proof turned 
to be characteristic--free, a somehow unusual fact in singularity theory.

The main result thus generalizes previous work of Zariski 
(\cite{Z}), Hironaka (\cite{H}), Abhyankar (\cite{GP}) and others, who 
treated specific cases tailored for their purposes on resolution of 
surface singularities. Note that, in fact, all these results were previous 
to Mulay's which is a key stone in our strategy. As the geometric behaviour of 
the equimultiple locus of three-dimensional varieties remains unknown (see,
for instance, \cite{Sp}), we hope this result can be used as a 
first step to this much more difficult and intriguing problem.

\section{Notation and technical results}

For the sake of completeness (and for the convenience of the
reader), we recall the basic facts and technical results related
to quadratic and monoidal transformations that were used in
\cite{PT} and which will be used afterwards.

For all what follows, let $\cS$ be an embedded algebroid surface
of multiplicity $n$,
$$
F = Z^n + \sum_{k=0}^{n-1} \left( \sum_{i,j} a_{ijk}X^iY^j \right)
Z^k = Z^n + \sum_{k=0}^{n-1} a_k(X,Y) Z^k
$$
a Weierstrass equation of $\cS$. We will note
$$
N_{\{X,Y,Z\}}(F) = \left\{ (i,j,k) \in \NN^3 \; | \; a_{ijk} \neq 0
\right\},
$$
\noindent although we will omit the subscript whenever the variables
is clear from the context.

\df
The elements of $\cE (\cS)$ different from $M$ will be called
equimultiple curves. The elements of $\cE_0(\cS)$ other than $M$ will be called
permitted curves.

\obs In particular, any $P \in \cE_0 (\cS)$ can be assumed to be,
for instance $(Z,X)$, after a suitable change of variables.
Clearly $(Z,X) \in \cE_0 (\cS)$ is equivalent to  $i+k \geq n$ 
for all $(i,j,k)\in N(F)$.

The monoidal transform of $\cS$,
centered in $(X,Z)$, in the point corresponding to the direction
$(\alpha:0 :\gamma)$ (say $\alpha \neq 0$) of the exceptional
divisor is the surface $\cS^{(1)}$ defined by the equation
$$
F^{(1)} = \left( Z_1 + \frac{\gamma}{\alpha} \right)^n + 
\sum_{(i,j,k) \in N(F)} a_{ijk}X_1^{i+k-n} Y_1^j
\left( Z_1 + \frac{\gamma}{\alpha} \right)^k.
$$

Observe that this
only makes sense (that is, gives a non--unit) whenever
$\ol{F} (\alpha,0,\gamma)=0$. The homomorphism
\begin{eqnarray*}
\pi^P_{(\ol{\alpha}:0:\gamma)}: K[[X,Y,Z]] & \longrightarrow &
K[[X_1,Y_1,Z_1]] \\
X & \longmapsto & X_1 \\
Y & \longmapsto & Y_1 \\
Z & \longmapsto & \displaystyle X_1 \left( Z_1 + \frac{\gamma}{\alpha}
\right) \\
\end{eqnarray*}
will be called the homomorphism associated to the monoidal
transformation in $(\ol{\alpha}:0:\gamma)$ or, in short,
the equations of the monoidal transformation. The overline is
because one must privilege a non-zero coordinate, but all the
possibilities define associated equations. 

The quadratic transform (that is, blowing--ups with center $M$)
in the point corresponding to the direction
$(\alpha:\beta:\gamma)$ (say $\alpha \neq 0$) of the exceptional
divisor is the surface $\cS^{(1)}$ defined by the equation
$$
F^{(1)} = \left( Z_1 + \frac{\gamma}{\alpha} \right)^n + 
\sum_{(i,j,k) \in N(F)} a_{ijk}X_1^{i+j+k-n} \left( Y_1
+ \frac{\beta}{\alpha} \right)^j \left( Z_1 + \frac{\gamma}{\alpha}
\right)^k.
$$

Again this only makes sense whenever
$\ol{F} (\alpha, \beta, \gamma)=0$. Analogously, the homomorphism
\begin{eqnarray*}
\pi^M_{(\ol{\alpha}:\beta:\gamma)}: K[[X,Y,Z]] & \longrightarrow &
K[[X_1,Y_1,Z_1]] \\
X & \longmapsto & X_1 \\
Y & \longmapsto & \displaystyle X_1 \left( Y_1 + \frac{\beta}{\alpha}
\right) \\
Z & \longmapsto & \displaystyle X_1 \left( Z_1 + \frac{\gamma}{\alpha}
\right) \\
\end{eqnarray*}
will be called the homomorphism associated to the quadratic
transformation in $(\ol{\alpha}:\beta:\gamma)$ or the equations of
the quadratic transformation. 

\obs
In the previous situation, consider a change of
variables in $K[[X,Y,Z]]$ given by
$$
\left\{ \begin{array}{ccc}
\varphi (X) & = & a_1 X' + a_2 Y' + a_3 Z' + \varphi_1 (X',Y',Z') \\
\varphi (Y) & = & b_1 X' + b_2 Y' + b_3 Z' + \varphi_2 (X',Y',Z') \\
\varphi (Z) & = & c_1 X' + c_2 Y' + c_3 Z' + \varphi_3 (X',Y',Z') \\
\end{array} \right.,
$$
with $\ord \left( \varphi_i \right) \geq 2$.

Assume also that
$$
\left\{ \begin{array}{ccc}
\alpha & = & a_1 \alpha' + a_2 \beta' + a_3 \gamma' \\
\beta & = & b_1 \alpha' + b_2 \beta' + b_3 \gamma' \\
\gamma & = & c_1 \alpha' + c_2 \beta' + c_3 \gamma' \\
\end{array} \right.
$$
with, say, $\gamma'\neq 0$. Then there is a unique
change of variables $\psi: K[[X_1,Y_1,Z_1]] \longrightarrow
K[[X'_1,Y'_1,Z'_1]]$ such that
$$
\psi \pi^M_{(\ol{\alpha}:\beta:\gamma)} = \pi^M_{(\alpha':\beta':
\ol{\gamma'})} \varphi.
$$
 
\df
Let $Q \in \cE(\cS)$, with $Q=(Z+H(X,Y),G(X,Y))$. Then for 
$\underline{u} \in \PP^2(K)$, the ideal
$$
\varpi^M_{\underline{u}} (Q) = \left( \frac{\pi^M_{\underline{u}}
(Z+H(X,Y))}{X_1}, \frac{\pi^M_{\underline{u}}(G(X,Y))}{X_1^{\ord(G)}} 
\right)
$$
is called the (strict) quadratic transform of $Q$ in the point
$\underline{u}$.

Obviously, this definition makes sense only if the quadratic
transform in the direction $\underline{u}$ does. There is a
natural version of monoidal transform of $Q$ with center $P$,
for all $P \in \cE_0(\cS)$.

\ntt
We will note by $\nu$ the natural isomorphism
$$
\nu : K[[X,Y,Z]]  \longrightarrow  K \left[\left[ X_1,Y_1,Z_1
\right]\right]
$$
sending $X$ to $X_1$, $Y$ to $Y_1$ and $Z$ to $Z_1$.

\section{The theorem}

We will restrict ourselves to the case which is interesting for
desingularization issues: that where $\cS$ and $\cS^{(1)}$ have
the same multiplicity. This leaves out some situations.

\lema If the tangent cone of $\cS$ is not a plane, the multiplicity of any
monoidal transform is strictly less than $n$.

\dem See \cite{PT} for a characteristic--free proof.

\obs In particular, note that if the tangent cone is not a plane,
there cannot be more than one permitted curve. In fact, assume we
have two curves, $P$ and $Q$, and choose $Z$ to be a regular
parameter with $P=(Z,G(X,Y))$, $Q=(Z,H(X,Y))$. If $F$ is a
Weierstrass equation with the usual form, then $P$ is permitted if
and only if $G^{n-k}|a_k$ for all $k=0,...,n-1$. As the same goes
for $Q$, it is clearly impossible that there exists some $a_k$
with $\ord(a_k)=n-k$. Hence the tangent cone must be a plane (in
fact it must be $Z=0$).

\thr Let $\cS$ be an algebroid surface and $\cS^{(1)}$ a quadratic or monoidal 
transform of $\cS$ having the same multiplicity.

\noindent
{\bf (a)} Let $\cS^{(1)}$ be the monoidal transform of $\cS$ with
center $P \in \cE_0 (\cS)$ then, either $\cE_0
\left( \cS^{(1)} \right) = \nu \left( \cE_0 (\cS) \right)$ or $\cE_0
\left( \cS^{(1)} \right) = \nu \left( \cE_0 (\cS) \setminus \{ P \}
\right)$.

\noindent
{\bf (b)} Let $\cS^{(1)}$ be the quadratic transform of $\cS$
in the point $\underline{u}$. 

{\bf (b.1)} If the tangent cone is not a plane then $\cE_0
\left( \cS^{(1)} \right)= \varpi^M_{\underline{u}} \left( \cE_0 (\cS) \right)$.

{\bf (b.2)} If the tangent cone is a plane then we can find 
three types of curves in $\cE_0 \left( \cS^{(1)} \right)$:

(i)  The exceptional divisor of the transform.

(ii) Primes $\varpi^M_{\underline{u}} (Q)$, with $Q \in \cE (\cS) \setminus
\cE_0(\cS)$, which are tangent to the exceptional divisor.

(iii) Primes $\varpi_{\underline{u}}^M (Q)$, with $Q \in \cE_0 (\cS)$, where 
both $\nu(Q)$ and $\varpi^M_{\underline{u}} (Q)$ are transversal to the 
exceptional divisor.

Moreover, if $E_0 \left( \cS^{(1)} \right)$ contains a prime of type (ii), then it 
also contains a prime of type (i).

\dem Although some partial results are common to the
characteristic $0$ case, we will repeat them or, at the very
least, we will give a detailed outline when appropriate for the
convenience of the reader. In what follows let $F$ be, as usual, a 
Weierstrass equation of $\cS$.

\vs

\noindent \fbox{{\bf Case (a)}}

\vs

This case presents the first (small) differences of argumentation
with the characteristic $0$, as the reader can check with
\cite{PT}. We have to prove that, after a monoidal transformation
with center $P \in \cE_0 (\cS)$, $Q \in \cE_0 (\cS)$ if and only if
$\nu (Q) \in \cE_0 \left( \cS^{(1)} \right)$, except maybe for $Q=P$.

Hence assume $Z$ is a parameter verifying that for all $I \in \cE
(\cS)$, $Z \in I$. After a change of variables in $K[[X,Y]]$, we
can assume $P$ to be $(Z,X)$ and $Q$ (other than $P$) to be $(Z,G(X,Y))$ with
$\ord (G)=1$. As we noticed above $\ol{F} = Z^n$, hence there is only
one direction in the exceptional divisor, $(1:0:0)$. An equation
for $\cS^{(1)}$ is then
$$
F^{(1)} = Z_1^n + \sum_{k=0}^{n-1} \frac{a_k(X_1,Y_1)}{X_1^{n-k}}
Z_1^k = Z_1^n + \sum_{k=0}^{n-1} a^{(1)}_k(X_1,Y_1) Z_1^k,
$$
and therefore $G(X_1,Y_1)^{n-k}|a_k^{(1)}(X_1,Y_1)$ if and only if
$G(X,Y)^{n-k}|a_k(X,Y)$. This finishes the case.

\vs

\noindent \fbox{{\bf Case (b.2)}}

\vs

Much like in zero characteristic, this case gives the basis for the
other one. It is also the point where the main differences
between both cases become notorious.

\obs First of all note that we can restrict ourselves to the case
where the tangent cone is $Z=0$ (this is as before an easy change
of variables) and the direction in the exceptional divisor is
$(1:0:0)$. If this is not the case, for the results at points
$(1:\alpha:0)$ it suffices considering the (commutative) diagram

\unitlength=.9cm
\begin{picture}(11.5,5)(.5,-1)

\put(2,3){\makebox(0,0){$K[[X,Y,Z]]$}}
\put(10,3){\makebox(0,0){$K[[X',Y',Z']]$}}
\put(2,0){\makebox(0,0){$K[[X_1,Y_1,Z_1]]$}}
\put(10,0){\makebox(0,0){$K[[X'_1,Y'_1,Z'_1]]$}}

\put(3.5,3){\vector(1,0){5}}
\put(3.5,0){\vector(1,0){5}}
\put(2.5,2.5){\vector(0,-1){2}}
\put(9.5,2.5){\vector(0,-1){2}}

\put(6,3.25){\makebox(0,0){$\scriptstyle \varphi$}}
\put(6,0.25){\makebox(0,0){$\scriptstyle \psi$}}
\put(2.25,1.5){\makebox(0,0)[r]{$\scriptstyle \pi^M_{(1:\alpha:0)}$}}
\put(9.75,1.5){\makebox(0,0)[l]{$\scriptstyle \pi^M_{(1:0:0)}$}}

\end{picture}

\noindent
with $\varphi$ given by
$$
\left\{ \begin{array}{ccl}
\varphi (X) & = & X' \\
\varphi (Y) & = & Y' + \alpha X' \\
\varphi (Z) & = & Z'
\end {array} \right.
$$

Of course the results at $(0:1:0)$ are clearly symmetric.

\obs Assume that we have a permitted curve in $\cS^ {(1)}$ of the
general type, say, $Q = \left( G_1, G_2 \right)$, with
$$
\begin{array}{clcl}
G_1 & = & \alpha_1 X_1 + \beta_1 Y_1 + \gamma_1 Z_1 + G'_1 \left( X_1,Y_1,Z_1
\right),& \mbox{ where } \ord \left( G'_1 \right)> 1 \\
G_2 & = & \alpha_2 X_1 + \beta_2 Y_1 + \gamma_2 Z_1 + G'_2 \left( X_1,Y_1,Z_1
\right), & \mbox{ where } \ord \left( G'_2 \right)> 1 \\
\end{array}
$$

As the multiplicity remains the same, the monomial $Z_1^n$ must
appear in $\ol{F^{(1)}}$. Hence either $\gamma_1$ or $\gamma_2$ must be non
zero. Let us suppose it is $\gamma_1 \neq 0$ and so we can substitute
$G_1$ by its associated Weierstrass polynomial with respect to
$Z_1$, of the form $Z_1+ a(X_1,Y_1)$ with $\ord(a) \geq 1$.
Now we change $Z_1$ by $-a(X_1,Y_1)$ in $G_2$ to obtain
$$
Q = \left( Z_1 + a \left( X_1,Y_1 \right),
\alpha X_1 + \beta Y_1 + b \left( X_1,Y_1 \right) \right),
$$
with $\ord(a), \ord(b)>1$.

We will look first at permitted curves in $\cS^{(1)}$ which are
transversal to the exceptional divisor.

\lema Under the hypothesis of case (b.2) assume that there is a
permitted curve $Q \in \cE_0 \left( \cS^{(1)} \right)$ which is
transversal to the exceptional divisor. Then there exists some $P
\in \cE_0 (\cS)$ such that $Q = \varpi_{(1:0:0)}^M (P)$ and $\nu (P)$ is
transversal to the exceptional divisor.

\dem In the above notation, $Q$ is transversal whenever $\beta \neq
0$. In this case we can change $\alpha X_1 + \beta Y_1 +
b(X_1,Y_1)$ for its associated Weierstrass polynomial with respect
to $Y_1$ and then make the corresponding substitution in
$a(X_1,Y_1)$ to obtain
$$
Q = \left( Z_1 + a' \left( X_1 \right),
Y_1 + b' \left( X_1 \right) \right),
$$
with $\ord \left(a' \right)>1$, $\ord \left( b' \right)>0$. Consider then 
the following diagram
\begin{picture}(11.5,5)(.5,-1)

\put(2,3){\makebox(0,0){$K[[X,Y,Z]]$}}
\put(10,3){\makebox(0,0){$K[[X',Y',Z']]$}}
\put(2,0){\makebox(0,0){$K[[X_1,Y_1,Z_1]]$}}
\put(10,0){\makebox(0,0){$K[[X'_1,Y'_1,Z'_1]]$}}

\put(3.5,3){\vector(1,0){5}}
\put(3.5,0){\vector(1,0){5}}
\put(2.5,2.5){\vector(0,-1){2}}
\put(9.5,2.5){\vector(0,-1){2}}

\put(6,3.25){\makebox(0,0){$\scriptstyle \varphi$}}
\put(6,0.25){\makebox(0,0){$\scriptstyle \psi$}}
\put(2.25,1.5){\makebox(0,0)[r]{$\scriptstyle \pi^M_{(1:0:0)}$}}
\put(9.75,1.5){\makebox(0,0)[l]{$\scriptstyle \pi^M_{(1:0:0)}$}}

\end{picture}

\noindent
with changes of variables
$$
\left\{ \begin{array}{ccl}
\varphi (X) & = & X' \\
\varphi (Y) & = & Y' - X'b' \left( X' \right)\\
\varphi (Z) & = & Z' - X'a' \left( X' \right)
\end {array} \right. , \quad
\left\{ \begin{array}{ccl}
\psi \left( X_1 \right)& = & X'_1 \\
\psi \left( Y_1 \right) & = & Y'_1 - b' \left( X'_1 \right)\\
\psi \left( Z_1 \right) & = & Z'_1 - a' \left( X'_1 \right)
\end {array} \right.
$$

As $\psi(Q)= \left( Z'_1,Y'_1 \right)$, we know that $\left(
Z'_1,Y'_1 \right)$ is permitted in $\cS^{(1)}$. But, looking at
the equations of the transformation $\pi^M_{(1:0:0)}$ this clearly
implies that $(Z',Y')$ was permitted in $\cS$. Therefore
$$
P = \varphi^{-1} \left(Z',Y' \right)= \left( Z+Xa'(X),Y+Xb'(X) \right)
$$
was permitted in $\cS$. It is clear that $\nu(P)$ is transversal to
the exceptional divisor and $\varpi_{(1:0:0)}^M (Q) = P$. This proves 
the lemma.

Now we will prove that the existence of permitted curves tangent
to the exceptional divisor implies that the exceptional divisor
lies in $\cE_0 \left( \cS^{(1)} \right)$. We will do it without
using the fact that $\ol{F}$ is the power of a linear form, and so
it will be still valid for case (b.1). 

\lema Under the hypothesis of case (b) assume that there is a
permitted curve $P \in \cE_0 \left( \cS^{(1)} \right)$ which is
tangent to the exceptional divisor. Then the exceptional divisor 
lies in $\cE_0 \left( \cS^{(1)} \right)$.

\dem From the remarks at the beginning of this case we can already
assume that $P$ has the form
$$
P= \left( Z_1+ a(Y_1), X_1+ b (Y_1) \right), \mbox{ with }
\ord(a), \ord(b) \geq 2.
$$

The proof is considerably different, depending on whether $a(Y_1)$
is zero or not. Let us assume $a(Y_1)=0$ (this is the easy part)
and let us write $F^{(1)}$ as usual:
$$
F^{(1)} = Z_1^n + \sum_{k=0}^{n-1} a_k^{(1)}
\left( X_1,Y_1 \right)Z_1^k,
$$
with the following decomposition
$$
a_k^{(1)} \left(X_1,Y_1 \right)=  b_k \left( X_1,Y_1 \right)
\left( X_1+b(Y_1) \right)^{n-k}, \;k=0,...,n-1,
$$
\noindent where $b_k$ does not divide $X_1 + b(Y_1)$.

Let us take any $k \in \{0,...,n-1\}$, let $t = \ord (b) \geq 2$
and let us call $X_1^rY_1^s$ the smallest monomial with respect to
the lexicographic order appearing in $b_k \left( X_1,Y_1 \right)$. Then the
monomial $X_1^rY_1^{s+t(n-k)}Z_1^k$ must occur in $a_k^{(1)}
\left( X_1,Y_1 \right) Z_1^k$.

Now, as this monomial appears after a quadratic transform in the point
$(1:0:0)$ of the exceptional divisor it is clear that it must hold
$$
r \geq s+t(n-k)+k-n \geq s+2(n-k)-(n-k) \geq n-k,
$$
and therefore $X_1^{n-k} | b_k \left( X_1,Y_1 \right)$ for all
$k$. This proves that the exceptional divisor is also permitted.

Let us move now to the more complicated case, when $a \left(Y_1
\right)\neq 0$. Now we write $F^{(1)}$ in the following form
$$
F^{(1)} = \left( Z_1 + a(Y_1) \right)^n+ \sum_{k=0}^{n-1} b_k
\left( X_1,Y_1 \right) \left( X_1 + b(Y_1) \right)^{n-k} \left(
Z_1 + a(Y_1) \right)^k,
$$
for some $b_k \left( X_1,Y_1 \right) \in K[[X_1,Y_1]]$. Hence the 
summand $a \left( Y_1 \right)^n$, which is a power series in 
$K[[Y_1]]$ of order strictly greater than $n$, appears in
the independent term.

But as $F^{(1)}$ comes from $F$ after a quadratic transformation
at the point $(1:0:0)$ of the exceptional divisor, there can be no
monomials in $K[[Y_1]]$ of order strictly greater then $n$ in the 
independent term, which implies
clearly that $b \left( Y_1 \right)$ cannot be $0$.

Let us write $a \left( Y_1 \right) = Y_1^r u \left( Y_1 \right)$,
$b \left( Y_1 \right) = Y_1^s v \left(Y_1 \right)$, for some $u,v$
units in $K[[Y_1]]$. Then, as previously, it must hold
$$
Y_1^{rn} u \left( Y_1 \right)^n + \sum_{k=0}^{n-1} b_k \left( 
0,Y_1 \right) Y_1^{rk} u \left( Y_1 \right)^k Y_1^{s(n-k)} 
v \left( Y_1 \right)^{n-k} = 0,
$$
which in particular implies $rn \geq sn$, that is, $r \geq s$.

Let us write $r = r's+t$ with $t \in \{0,...,s-1\}$ and $r' \geq 1$. 
Then we can rewrite our ideal as
$$
\begin{array}{ccl}
P & = & \displaystyle \left( Z_1 \pm Y_1^tX_1^{r'} \left( \frac{u(Y_1)}
{v(Y_1)^{r'}} \right), X_1 + Y_1^s v(Y_1) \right) \\ \\
& = & \left( Z_1 + X_1^{r'} c(Y_1), X_1+Y_1^sv(Y_1) \right),
\end{array}
$$
or, alternatively
$$
P = \left( Z_1 + X_1 g \left( X_1,Y_1 \right), X_1 + b(Y_1) \right).
$$

With this new form it must hold
$$
\begin{array}{l}
F^{(1)} = \left( Z_1 + X_1 g(X_1,Y_1) \right)^n + \\ \\
\displaystyle \quad \quad \quad \quad + \sum_{k=0}^{n-1} c_k
\left( X_1,Y_1 \right) \left( X_1 + b(Y_1) \right)^{n-k} \left(
Z_1 + X_1 g(X_1,Y_1) \right)^k.
\end{array}
$$

We want to prove $X_1^{n-k}|c_k \left(X_1,Y_1\right)$ for $k=0,...,n-1$.
If this were not so, let $k_0$ be the biggest index for which this does
not happen. Let us write $X_1^pY_1^q$ the smaller monomial for the 
lexicographic ordering appearing in $c_{k_0} \left( X_1,Y_1 \right)$.
Note that $p<n-k_0$. Then the monomials in
$$
X_1^p Y_1^q b \left( Y_1 \right)^{n-k_0} Z_1^{k_0}
$$ 
cannot cancel with any others coming from the expansions of the remaining 
summands. This is so because, if $k>k_0$ the monomials in the $k$--th 
summand are in $\left( Z_1,X_1 \right)^n$ and if $k<k_0$ the monomials 
in the $k$--th summand have smaller exponent in $Z_1$. Hence the monomial
$$
X_1^p Y_1^{q+s(n-k_0)}Z_1^{k_0}
$$
actually appears in $F^{(1)}$. As $F^{(1)}$ is the result of a quadratic
transform in $(1:0:0)$, we have
$$
p \geq q + s \left( n-k_0 \right) + k_0 - n = q + (s-1) \left( n- k_0 
\right)  \geq q + n - k_0 \geq n - k_0,
$$
which contradicts our assumption. This finishes the proof of the lemma.

\vs

Note that, from the previous arguments, a curve which is tangent to the
exceptional divisor can also be written as
$$
P = \left( Z_1 + h \left( X_1, Y_1 \right), X_1 + b \left( Y_1 \right)
\right),
$$
where $b \left( Y_1 \right) \in K[[Y_1]]$ and $h \left( X_1,Y_1 \right)
\in K[[X_1]][Y_1]$, with $\ord (b) = s \geq 2$ and $\deg_{Y_1}(h) < s$.
This will be our preferred form in what follows.

In order to prove the existence of a singular equimultiple curve
in $\cE(\cS)$ which is taken into the permitted curve that is
tangent to the exceptional divisor, we will begin by proving that
we can find a convenient parameter to work with in this
situation.

\lema In the situation of case (b.2), let us assume that, after
the quadratic transform in the direction $(1:0:0)$ of the
exceptional divisor, we have the following permitted curve in
$\cS^{(1)}$:
$$
P = \left( Z_1 + h \left(X_1,Y_1 \right),
X_1 + b(Y_1) \right),
$$
with $\ord(b) = s>1$ and $h \left(X_1,Y_1 \right) \in K [[X_1]][Y_1]$ 
with $\deg_{Y_1} \left( h \right) < s$.

Then there exists an isomorphism $\varphi: K[[X,Y,Z]]
\longrightarrow K[[X',Y',Z']]$ which makes the following
diagram commutative

\unitlength=.9cm
\begin{picture}(11.5,5)(.5,-1)

\put(2,3){\makebox(0,0){$K[[X,Y,Z]]$}}
\put(10,3){\makebox(0,0){$K[[X',Y',Z']]$}}
\put(2,0){\makebox(0,0){$K[[X_1,Y_1,Z_1]]$}}
\put(10,0){\makebox(0,0){$K[[X'_1,Y'_1,Z'_1]]$}}

\put(3.5,3){\vector(1,0){5}}
\put(3.5,0){\vector(1,0){5}}
\put(2.5,2.5){\vector(0,-1){2}}
\put(9.5,2.5){\vector(0,-1){2}}

\put(6,3.25){\makebox(0,0){$\scriptstyle \varphi$}}
\put(6,0.25){\makebox(0,0){$\scriptstyle \psi$}}
\put(2.25,1.5){\makebox(0,0)[r]{$\scriptstyle \pi^M_{(1:0:0)}$}}
\put(9.75,1.5){\makebox(0,0)[l]{$\scriptstyle \pi^M_{(1:0:0)}$}}

\end{picture}

\noindent
with $\psi$ given by
$$
\left\{ \begin{array}{ccl}
\psi \left( X_1 \right) & = & X_1' \\
\psi \left( Y_1 \right) & = & Y_1' \\
\psi \left( Z_1 \right) & = & \displaystyle
Z_1' - h \left( X'_1, Y'_1 \right)
\end {array} \right.
$$

\dem Let us write
$$
h \left( X_1,Y_1 \right) = h_1 \left( X_1,Y_1 \right) + h_2 \left(
X_1, Y_1 \right),
$$
where the monomials of $h_1$ are precisely those $X_1^\alpha Y_1^\beta$ 
with $1+\alpha \geq \beta$, all the remaining ones being in $h_2$. Note
that $h_2 \in K[X_1,Y_1]$.

If $h_2 = 0$ then we get the result we look for with
$$
\left\{ \begin{array}{ccl}
\varphi \left( X \right) & = & X' \\
\varphi \left( Y \right) & = & Y' \\
\varphi \left( Z \right) & = & \displaystyle
Z' - X' h_1 \left( X', \frac{Y'}{X'} \right)
\end {array} \right.
$$

If $h_2 \neq 0$ then the previous isomorphism $\varphi$, together with
$$
\left\{ \begin{array}{ccl}
\psi \left( X_1 \right) & = & X_1' \\
\psi \left( Y_1 \right) & = & Y_1' \\
\psi \left( Z_1 \right) & = & \displaystyle
Z_1' - h_1 \left( X'_1,Y'_1 \right)
\end {array} \right.
$$
makes the diagram commute. In this case, we would have the permitted curve
$$
\psi (P) = \left( Z'_1 + h_2 \left(X'_1,Y'_1 \right), 
X'_1 + b (Y'_1) \right)
$$
in the equimultiple locus of the surface ${\cal S}^{(1)}$. The equation for 
${\cal S}^{(1)}$, in $\left\{ X_1',Y_1',Z_1' \right\}$ can then be written as
$$
\begin{array}{l}
F^{(1)} =  \left( Z_1' + h_2 \left( X'_1,Y'_1 \right) \right)^n + \\ \\
\displaystyle \quad \quad \quad \quad \quad + \sum_{k=0}^{n-1} b_k 
\left( X'_1,Y'_1 \right) \left( X'_1 + b (Y'_1) \right)^{n-k} 
\left( Z'_1 + h_2 \left( X'_1,Y'_1 \right) \right)^k.
\end{array}
$$

For the sake of simplicity, we will rename our variables as $\{X,Y,Z\}$, as 
no further changes of variables are needed from this point onwards.

Assume $(j,d)$ is the smallest exponent in $h_2(X,Y)$ for the lexicographic 
ordering verifying $1+j < d \leq s-1$. Let us prove
$$
X^{d+(n-k-1)j}|b_k(X,Y), \mbox{ for all } k=0,...,n-1. 
$$

\vs

{\sc Step 1.--} Let us begin by the case $k=n-1$; we want to prove 
$X^d|b_{n-1}(X,Y)$.

\vs

The coeffcient of $Z^{n-1}$ in $F^{(1)}$ must be
$$
\left( \begin{array}{c} n \\ 1 \end{array} \right) h_2 + b_{n-1}(X,Y)(X+b(Y)).
$$

Let $(\alpha, \beta)$ be the minimal exponent appearing in $b_{n-1}(X,Y)$
with respect to the lexicographic ordering. Then $(\alpha, \beta+s)$ is
the minimal exponent appering in $b_{n-1}(X,Y)(X+b(Y))$ and it cannot cancel
with any monomial in $h_2$ as $\deg_{Y}\left( h_2 \right) \leq s-1$. Hence
$(\alpha, \beta+s,n-1) \in N \left( F^{(1)} \right)$ and, since $F^{(1)}$ is
the equation for a quadratic transform in $(1:0:0)$, it must hold
$$
\alpha \geq \beta + s+ (n-1) - n = \beta + s-1 \geq s-1 \geq d.
$$

\vs

{\sc Step 2.--} Let us assume the result is true for $k, k+1,...,n-1$ and
let us prove 
$$
X^{d+(n-k)j}|b_{k-1} (X,Y).
$$

We fix now our attention in the coefficient of $Z^{k-1}$ in $F^{(1)}$, which is
formed by the expansion of 
$$
b_{k-1}(X,Y) (X+b(Y))^{n-k+1},
$$
and also by monomials coming from the expression
$$
\left( Z+h_2 \right)^n + (X+b) b_{n-1} \left( Z+h_2 \right)^{n-1}+...+
(X+b)^{n-k} b_k \left( Z + h_2 \right)^k.
$$

If we mimic the case above and denote by $(\alpha, \beta)$ the minimal exponent 
appearing in $b_{k-1} (X,Y)$ for the lexicographic ordering, then the minimal
exponent for $b_{k-1}(X+b)^{n-k+1}$ must be $(\alpha, \beta+(n-k+1)s)$. And,
furthermore, this exponent cannot be cancelled with any other from the 
developement of the summands of the expression above. Let us see this with closer
detail.

Clearly, it cannot cancel with any monomial from 
$$
\left( \begin{array}{c} n \\ k-1 \end{array} \right) h_2 (X,Y)^{n-k+1},
$$
as all of their degrees with respect to $Y$ are smaller than $(s-1)(n-k+1)$.

Now, if $\alpha < d+(n-k)j$ there cannot be cancellation with monomials of
$$
\left( \begin{array}{c} n-m \\ k-1 \end{array} \right)b_{n-m} (X+b)^m h_2^{n-m-k+1},
\mbox{ (for $m=1,...,n-k$)},
$$
as their orders with repect to $X$ are greater or equal (by induction
hypothesis) than 
$$
d+(m-1)j+j(n-m-k+1)=d+j(n-k).
$$

Hence $(\alpha, \beta+(n-k+1)s, k-1) \in N \left( F^{(1)} \right)$ and, as before,
it must then hold
$$
\begin{array}{lclcl}
\alpha & \geq & \beta+(n-k+1)s+k-1-n & \geq &(n-k+1)(s-1) \\
& \geq & (n-k+1)d & = & d+(n-k)d \\
& > & d+(n-k)j,
\end{array}
$$
which is a contradiction. Hence $X^{d+(n-k-1)j}|b_k(X,Y)$, for all 
$k=0,...,n-1$. 

Now the independent term of $F^{(1)}$ is
$$
F^{(1)}(X,Y,0)= h_2^n + \sum_{k=0}^{n-1} b_k(X+b)^{n-k}h_2^k.
$$

The first summand features the exponent $(nj,nd,0)$, which must be in 
$N \left( F^{(1)} \right)$, as the order with respect to $X$ of any
other monomial is, for some $k\in \{0,...,n-1\}$ greater or equal than
$$
d+(n+k-1)j+kj = d+(n-i)j > nj+1 > nj.
$$

However, this monomial cannot appear in $F^{(1)}$, as $j<d$ and $F^{(1)}$
is the equation of a quadratic transform in $(1:0:0)$. This shows $h_2=0$ and
finishes the proof of the lemma.

\obs After this lemma we can assume that, in our situation, $Z_1$
lies in both the exceptional divisor and $P$. So, for finishing our
case, we are given the surface $\cS$ defined by
$$
F = Z^n +\sum_{k=0}^{n-1} a_k (X,Y) Z^k,
$$
and its quadratic transform $\cS^{(1)}$ on the direction $(1:0:0)$,
defined by
$$
F^{(1)} = Z_1^n +\sum_{k=0}^{n-1} a^{(1)}_k (X_1,Y_1) Z_1^k,
$$
about which we know the following facts:
\begin{itemize}
\item The multiplicity of $\cS$ and $\cS^{(1)}$ is the same.
\item There is a curve $P = (Z_1,X_1+G(Y_1)) \in \cE_0 \left(
\cS^{(1)} \right)$, with $\ord(G) \geq 2$.
\end{itemize}

We have to prove that there is a curve $Q \in \cE (\cS) \setminus
\cE_0 \left( \cS \right)$ such that $\varpi_{(1:0:0)}^M (Q) = P$. This was
proved in our previous paper \cite{PT} in a characteristic--free
way. In fact, the really difficult part in positive characteristic
is proving that we can assume $P$ to have this particular form,
and this has been done in the previous lemma.

We give a brief outline of the proof: first one shows that it is 
enough to prove that one can find a power series 
$H \left( X_1,Y_1 \right)$ verifying:
\begin{itemize}
\item $\ord (H) = \ord (G) = \lambda > 2$.
\item $H$ is regular in $Y_1$ of order $\lambda$.
\item There is a unit $u \left( X_1,Y_1 \right)$ such that
$$
\frac{1}{X_1^\lambda} H \left( X_1,X_1Y_1 \right) = u \left(
X_1,Y_1 \right) \left( X_1+ G (Y_1) \right).
$$
\end{itemize}

Then one proves that this power series actually exists in
a constructive way (in an ample sense, of course): it consists
simply in writing
$$
\begin{array}{lll}
X_1 + G(Y_1) & = & \displaystyle X_1 + \sum_{i \geq \lambda}
\alpha_i Y_1^i \\ \\
H (X_1,Y_1) & = & \displaystyle \sum_{k \geq \lambda} \left(
\sum_{i+j=k} \beta_{ij}X_1^iY_1^j \right) \\ \\
u (X_1,Y_1) & = & \displaystyle \sum_{k \geq 0} \left(
\sum_{i+j=k} \gamma_{ij}X_1^iY_1^j \right) \\ \\
\end{array}
$$
and then imposing the third condition above. In this way it is
straightforward seeing that one can, in fact, construct $H$
and $u$ with the desired properties.

\vs

\noindent \fbox{{\bf Case (b.1)}}

\vs

There are no great differences here with respect to the null
characteristic proofs. For instance, take $P =
(\alpha:\beta:\gamma)$ a direction in the exceptional divisor with
multiplicity $r$. Then the quadratic transform of $\cS$ on
$(\alpha:\beta:\gamma)$ has, at most, multiplicity $r$. This is
plain from the very definition of multiplicity.

So the hypothesis of (b.1) are filled only if the direction chosen
is one of multiplicity $n$. Eventually, changing the variables we
may consider that the point is $(0:1:0)$ (and $\ol{F} \in K[X,Z]$).

For proving that the quadratic transform cannot have new permitted
curves note that, in (b.2), we have shown that, if a new
permitted curve appears, so does the exceptional divisor (whether
$\ol{F}$ is the power of a linear form or not). But $(Z,Y)$ cannot
be a permitted curve, since $\ol{F^{(1)}}$ contains monomials in $K[X,Z]$
other than $Z^n$.

It is also clear that the quadratic transform does not erase
permitted curves either. If there is a permitted curve we may take
it to be $(Z,X)$, after a change of variables which does not affect 
$(0:1:0)$. Clearly this curve cannot disappear from the equimultiple 
locus after a quadratic transform on $(0:1:0)$.

This finishes the proof of the theorem.

\end{document}